\documentclass[reqno]{amsart}
\usepackage{amssymb,amsmath,amsthm,amsfonts,fancyhdr}
\usepackage{tikz}
\usepackage{indentfirst}
\usepackage{mathrsfs}
\usepackage{setspace}
\usepackage{color}
\usepackage{cite}
\usepackage{bm}
\numberwithin{equation}{section}

\newtheorem{thm}{Theorem}[section]

\newtheorem{lemma}[thm]{Lemma}
\newtheorem{cor}[thm]{Corollary}
\newtheorem{re}[thm]{Remark}
\newcommand{\osc}{{\mbox{osc}}}

  \newcommand{\dist}{\mbox{dist}}
  \newcommand{\diam}{{\mbox{diam}}}

  \newcommand{\supp}{{\mbox{supp}}}

  \numberwithin{equation}{section}
  \numberwithin{figure}{section}

\begin{document}

\title[${{W}^{2,p}}$ Estimates on ${C^{1,\alpha}}$ Domains]{ $\bm{{W}^{2,p}}$ Estimates for Elliptic Equations on $\bm{C^{1,\alpha}}$ Domains}

\author{Dongsheng Li}
\author{Xuemei Li}
\author{Kai Zhang}

\address{School of Mathematics and Statistics, Xi'an Jiaotong University, Xi'an, P.R.China 710049.}
\address{School of Mathematics and Statistics, Xi'an Jiaotong University, Xi'an, P.R.China 710049.}
\address{School of Mathematical Sciences, Shanghai Jiao Tong University, Shanghai, P.R.China 200240.}

\email{lidsh@mail.xjtu.edu.cn}
\email{pwflxm@stu.xjtu.edu.cn}
\email{zhangkaizfz@gmail.com}

\begin{abstract}
In this paper,  a new method is represented to investigate boundary $W^{2,p}$ estimates for elliptic equations, which is, roughly speaking, to derive boundary $W^{2,p}$ estimates from interior $W^{2,p}$ estimates by Whitney decomposition. Using it, $W^{2,p}$ estimates on $C^{1,\alpha}$ domains are obtained for nondivergence form linear elliptic equations and further more, fully nonlinear elliptic equations are also considered.
\end{abstract}



\maketitle

\section{Introduction}
\label{}

In this paper, we represent a new method of investigating boundary $W^{2,p}$ estimates for nondivergence form linear elliptic equations.
By virtue of Whitney decomposition, we derive local boundary $W^{2,p}$ estimates from interior $W^{2,p}$ estimates. This is a new point of view for boundary estimates and from it, we give a new proof of $W^{2,p}$ estimates on $C^{1,\alpha}$ domains as $1-1/p<\alpha\leq 1$ stated in \cite{KO}. It should be pointed out that our method can be extended to fully nonlinear elliptic equations as in Section 5.
Our result for linear equations is the following theorem.

\bigskip
\begin{thm}\label{t1.1}
Let $1<p<\infty$, $1-1/p<\alpha\leq 1$ and $\Omega$ be a bounded domain in $\mathbb{R}^n$ with a $C^{1,\alpha}$ boundary portion $T\subset \partial\Omega$. Let $L$ be an elliptic operator in nondivergence form:
\begin{equation*}
Lu=a^{ij}D_{ij}u+b^{i}D_{i}u+cu
\end{equation*}
with coefficients satisfying for some positive constants $0<\lambda\leq \Lambda<\infty$,
\begin{equation*}
\lambda|\xi|^2\leq a^{ij}(x)\xi_{i}\xi_{j}\leq\Lambda|\xi|^2 \ \ \mathrm{for} \ a.e.\  x\in \Omega \ \mathrm{and \ all} \ \xi\in \mathbb{R}^n
\end{equation*}
and
\begin{equation*}
\begin{array}{ccc}
 ||a^{ij}||_ {C(\bar \Omega)},\ ||b^i||_{L^\infty(\Omega)},\  ||c||_{L^\infty(\Omega)}\leq\Lambda.
\end{array}
\end{equation*}
Suppose that $u\in W^{2,p}(\Omega)$ is a strong solution of
$Lu=f$ in $\Omega$ and $u=g$ on $T$ in the sense of $W^{1,p}(\Omega)$ with $f\in L^{p}(\Omega)$ and $g\in W^{2,p}(\Omega)$. Then, for any domain $\Omega'\subset\subset\Omega\cup T$,
\begin{equation}\label{1.5}
||u ||_{W^{2,p}(\Omega')}\leq C\left( ||u ||_{L^{p}(\Omega)}+ ||f||_{L^{p}(\Omega)}+||g||_{W^{2,p}(\Omega)}\right),
\end{equation}
where $C$ depends on $n,\lambda,\Lambda,\alpha,p,T,\Omega',\Omega$ and the moduli of continuity of $a^{ij}$.
\end{thm}

\bigskip

\begin{re}\label{re}
$(i)$ If $T=\partial\Omega$ in Theorem \ref{t1.1}, we obtain a global $W^{2,p}(\Omega)$ estimate on the bounded $C^{1,\alpha}$ domain.

$(ii)$ The continuity assumption on $a_{ij}$ can be relaxed to $\osc_{B_r(x_0)\cap\Omega}a_{ij}\leq\delta$ for any small ball $B_r(x_0)$ and some suitably small $\delta>0$. By invoking the Sobolev embedding theorem, we may weaken the conditions on the lower order coefficients of $L$ to $b^i\in L^q(\Omega)$, $c\in L^r(\Omega)$, where $q>n$ if $p\leq n$, $q=p$ if $p>n$, $r>n/2$ if $p\leq n/2$, $r=p$ if $p>n/2$.

\end{re}

\bigskip

$W^{2,p}$ regularity plays a vital role in the regularity theory of partial differential equations. Interior $W^{2,p}$ estimates for Poisson's equation were first established in \cite{CZ} by explicit representation formulas involving singular integral operators.
Interior $W^{2,p}$ estimates for general nondivergence form elliptic equations are obtained on account of the fundamental observation that they can be treated locally as a perturbation of constant coefficient equations.
Later, Wang \cite{WL2} demonstrated a new proof via maximal function approach that is originated by Caffarelli \cite{C}.

 Boundary $W^{2,p}$ estimates are first established on flat domains by Schwarz reflection principle and then on $C^{1,1}$ domains alongside a flattening argument, where the $C^{1,1}$ regularity of domains is needed since the second order derivatives of the flattening mapping appear in the transformed elliptic operators (cf. \cite{GT}).

Based on theory of Sobolev multipliers, Maz'ya and Shaposhnikova \cite{MS} relaxed $C^{1,1}$ regularity of domains to
$M^{2-1/p}_{p}(\delta)$ (see Section 14.3.1 in \cite{MS} for its definition) for $1<p\leq n$ and $W^{2-1/p,p}$ for $n<p<\infty$, where $\delta$ depends on the moduli of continuity of $a^{ij}$.  They also proved as $1-1/p<\alpha\leq 1$, $C^{1,\alpha}\subset M^{2-1/p}_{p}(\delta)$ if
$1<p\leq n$ and $C^{1,\alpha}\subset W^{2-1/p}_{p}$ if
$n<p<\infty$. Moreover, they constructed $C^{1,1-1/p}$ domains where no solutions exist in $W^{2,p}$, which implies that the condition $1-1/p<\alpha$ can not be weakened.
For $p=2$, Kondrat'ev and $\grave{E}$idel'man [3] used the Fourier series technique to construct counterexamples showing that  $W^{2,2}$ estimates are invalid for $C^{1,1/2}$ domains. We refer to
Kondrat'ev and Oleinik \cite{KO} for a survey of the theory of boundary value problems in nonsmooth domains.

This paper investigates boundary $W^{2,p}$ estimates on $C^{1,\alpha}$ domains again not using singular integrals (and Sobolev multipliers).
The main idea is to derive local boundary $W^{2,p}$ estimates on $C^{1,\alpha}$ domains from interior $W^{2,p}$ estimates by Whitney decomposition.
Our approach is more direct and is applicable to both linear elliptic and fully nonlinear
elliptic equations. The proof is built upon Whitney decomposition, which is an effective tool for obtaining boundary estimates from interior estimates. For instance,  Cao, Li and Wang \cite{CL} utilized it to prove the optimal weighted $W^{2,p}$  estimates for elliptic equations with non-compatible conditions.

We illustrate our idea as follows. Let $\{Q_k\}_{k=1}^{\infty}$  be  Whitney decomposition of $\Omega_1$(Suppose $0\in\partial\Omega$ and denote $\Omega_r=\Omega\cap B_r$) and  $\widetilde Q_k=\frac65Q_k$ be $\frac65-$dilation of $Q_k$ with respect to its center. We suppose $L=\Delta$ and consider
\begin{equation*}
\left\{
\begin{array}{lr} \Delta u=f\ \ \ \ \ \ \mathrm{in}\ \ \ \ \ \ \Omega_1,\\ \ \ u=0\ \ \ \ \ \ \ \mathrm{on}\ \ \ (\partial\Omega)_1,
\end{array}
\right.
\end{equation*}
where $(\partial\Omega)_1=\partial\Omega\cap B_1.$ Deduce from interior $W^{2,p}$ estimates that
\begin{equation*}
\begin{aligned}
||D^2 u||_{L^{\tilde p}(Q_k)}^{\tilde p}
\leq C \left(d_k^{-2{\tilde p}}||u-l||_{L^{{\tilde p}}(\widetilde Q_k)}^{\tilde p}
+||f||_{L^{{\tilde p}}(\widetilde Q_k)}^{\tilde p}
\right)
\end{aligned}
\end{equation*}
for some ${\tilde p}>1$, some constant $C=C(n,{\tilde p})$ and any affine function $l$, where $d_k$ denotes the diameter of $Q_k$. If $C^{1,\tilde\alpha}$ estimate holds up to the boundary, we can take $l$ such that
$$|u(x)-l(x)|\leq  C\dist(x,\partial\Omega)^{1+\tilde\alpha},\ \forall x\in \Omega_{3/4}.$$
That is,
$$||u-l||^{\tilde p}_{L^{{\tilde p}}(\widetilde Q_k)}\leq Cd_k^{(1+\tilde\alpha){\tilde p}+n}.$$
It follows that
\begin{equation*}
\begin{aligned}
||D^2 u||_{L^{\tilde p}(Q_k)}^{\tilde p}
\leq C \left(d_k^{-{\tilde p}+\tilde\alpha {\tilde p}+n}
+||f||_{L^{{\tilde p}}(\widetilde Q_k)}^{\tilde p}
\right).
\end{aligned}
\end{equation*}
Take sum on both sides with respect to $k$ and if
$-{\tilde p}+\tilde\alpha {\tilde p}+1>0$, then $\sum_kd_k^{-{\tilde p}+\tilde\alpha {\tilde p}+n}$ is convergent (cf. Lemma 2.5).
Observe that here $\tilde\alpha=\min\{\alpha,1-\frac{n}{p}\}$ if we assume that $f\in L^{p}$ with $p>n$, and that ${\tilde p}<\frac{p}{n}$ since $0<-{\tilde p}+\tilde\alpha {\tilde p}+1\leq1-\frac{n\tilde p}{p}$.
Thus, we obtain a {\it rough} version $W^{2,p}$ estimate up to the boundary (cf. Remark \ref{hi} and Theorem \ref{t5.1}).
To obtain an {\it exquisite} version $W^{2,p}$ estimate up to the boundary (Theorem 1.1), we need further decompose $u$ and for linear equations, this is possible. Actually,
we set $u=v+w$ such that $v$ is a harmonic and
$w=\sum_{l}w_l$ with $w_l$ satisfying
\begin{equation*}
\left\{
\begin{array}{lr} \Delta w_l=f\chi_{Q_l}\ \ \mathrm{in}\ \ \ \ \Omega_1,\\ \ \ w_l=0\ \ \ \ \ \ \ \mathrm{on}\ \ \ \partial\Omega_1.
\end{array}
\right.
\end{equation*}
Since a large quantity of $w_l$ are harmonic in $Q_k$ and {\it better} boundary $C^{1,\alpha}$ estimates hold for them, we can improve the above {\it rough} estimate (cf. Section 4).

The paper is organized as follows. In Section 2, Whitney decomposition and its relevant properties are concluded.
 In Section 3, we demonstrate some basic estimates for elliptic equations including $W^{2,p}$ estimates for harmonic functions on $C^{1,\alpha}$ domains.
  In Section 4, we show local boundary $W^{2,p}$ estimates for Poisson's equation on $C^{1,\alpha}$ domains whose easy consequence is Theorem \ref{t1.1}.
  In Section 5, $W^{2,p}$ estimates for fully nonlinear elliptic equations on $C^{1,\alpha}$ domains are considered.

We end this section by listing some notations.
\bigskip

\noindent\textbf{Notation.} \\
1. $e_i=(0,...,0,1,...,0)=i^{th}$ standard coordinate vector.\\
2. $x'=(x^1,x^2,...,x^{n-1})$ and $x=(x',x^n).$\\
3. $\mathbb{R}^n_+=\{x\in \mathbb{R}^n:x^n>0\}.$\\
4. $B_r(x_0)=\{x\in \mathbb{R}^n: |x-x_0|<r\}$ and $B^+_r (x_0) = B_r(x_0) \cap  \mathbb{R}^n_+$.\\
5. $B_r'=\{x'\in \mathbb{R}^{n-1}: |x'|<r\}$ and $T_r= \{(x',0): x'\in B_r'\}.$\\
6. $\Omega_r (x_0)= \Omega\cap B_r(x_0)$ and $(\partial\Omega)_r (x_0)= \partial\Omega\cap B_r(x_0)$. We omit $x_0$ when $x_0=0.$
7. $\diam E$ = diameter of $E$,  $\forall E\subset\mathbb{R}^n.$\\
8. $\dist(E,F)$ = distance from $E$ to $F$, \ $\forall E,F\subset\mathbb{R}^n.$\\
9. $p'=p/(p-1)$ for $1<p<\infty$.

\section{Whitney decomposition}
In what follows, by a cube we mean a closed cube in $\mathbb{R}^n$, with sides parallel to the axes. We say two such cubes are disjoint if their interiors are disjoint.

\bigskip

\begin{lemma}\label{l2.1}(Whitney decomposition)
Let $F$ be a non-empty closed set in $\mathbb{R}^n$ and $\Omega$ be its complement. Then there exists a sequence of cubes $Q_k$(called the Whitney cubes of $\Omega$) such that

(i) $\Omega=\bigcup_{k=1}^{\infty}Q_{k}$;

(ii) The $Q_k$ are mutually disjoint;

(iii) $d_k\leq \mathrm{dist}(Q_k,F)\leq 4d_k$, where $d_k=\diam Q_k$.
\end{lemma}

\bigskip

\begin{lemma}\label{l2.2}
Let $Q_k$ be as in Lemma \ref{l2.1} and $\widetilde{Q}_k=\frac65{Q_k}$ be $\frac65-$dilation of $Q_k$ with respect to its center.
Then

$(i)$ $\Omega=\bigcup_{k=1}^{\infty}\widetilde Q_k$;

$(ii)$ Each point of $\Omega$ is contained in at most $12^n$ of the cubes $\widetilde Q_k$.
\end{lemma}

\bigskip

For the proof of the above two lemmas, we refer to Theorem 1 and Proposition 1-3 in Section VI.1 in \cite{ST}.

In the following, we always assume that $0\in \partial\Omega$ and there exists $\varphi\in C^{0,1}(B_{R}')$ with $||\varphi||_{C^{0,1}(B_{R}')}\leq K$ such that
\begin{equation*}
\Omega_R=\{x^n>\varphi(x')\}\cap B_R\ \ \mathrm{and}\ \ (\partial\Omega)_R=\{x^n=\varphi(x')\}\cap B_R
\end{equation*}
for some positive constants $R$ and $K$.
Throughout this paper, we assume $R=1$.
Let $\{Q_k\}_{k=1}^{\infty}$ be Whitney decomposition of $\Omega_{1}$ and  $\widetilde{Q}_k=\frac65 {Q_k}$.

\bigskip

\begin{lemma}\label{l2.3}
\begin{equation}\label{3.11}
\Omega_{r/3}\subset
\bigcup\limits_{\widetilde Q_k\subset \Omega_{r}}Q_k\ \ \mathrm{for}\ \ 0<r\leq1.
\end{equation}
\end{lemma}
\begin{proof}
If not, there exist a point $x\in \Omega_{r/3}$ and a cube $Q_k$ such that $x\in Q_k$ but $\widetilde Q_k\not\subset \Omega_{r}$. It follows that there exists a point  $y\in \widetilde Q_k$ with $|y|\geq r$. Then we deduce from Lemma \ref{l2.1} (iii) that
$$\dist (Q_k, \partial \Omega_{1})\geq \diam Q_k=\frac 56\diam \widetilde Q_k\geq \frac 56(|y|-|x|)\geq \frac 59 r.$$
Since $x\in Q_k\cap\Omega_{r/3}$ and $0\in \partial\Omega$,
$$\dist (Q_k, \partial \Omega_{1})\leq|x|\leq r/3.$$
Thus we get a contradiction.
\end{proof}

\bigskip
\begin{lemma}\label{d}
For any $x_0\in \bar\Omega_1$ and $r>0$ with $\Omega _{r}(x_0)\subset\Omega_{1}$,
\begin{equation}\label{bb}
|\Omega _{r}(x_0)\cap\{\dist(x,(\partial\Omega)_1)\leq d\}| \leq Cr^{n-1}d\ \ \mathrm{for}\ \ d>0.
\end{equation}
where $C$ depends only on $n$ and $K$.
\end{lemma}
\begin{proof}
Since $0\in \partial\Omega$, $\Omega_1=\{x^n>\varphi(x')\}\cap B_1$
and $(\partial\Omega)_1=\{x^n=\varphi(x')\}\cap B_1$
with $||\varphi||_{C^{0,1}(B_{R}')}\leq K$, we have
$$\begin{array}{lr}\Omega _{r}(x_0)\cap\{\dist(x,(\partial\Omega)_1)\leq d\}
\subset\{|x'-x_0'|\leq r,\varphi(x')\leq x^n\leq \varphi(x')+(K+1)d\}.
\end{array}$$
Since $|\{|x'-x_0'|\leq r,\varphi(x')\leq x^n\leq \varphi(x')+(K+1)d\}|\leq Cr^{n-1}d$, we have
$$|\Omega _{r}(x_0)\cap\{\dist(x,(\partial\Omega)_1)\leq d\}| \leq Cr^{n-1}d,$$
where $C$ depends only on $n$ and $K$.
\end{proof}

For further calculation, we set
\begin{equation}\label{fs}
\mathcal{F}^{s}=\bigcup_{k}\ \{Q_k:2^{-s-1}< d_k\leq 2^{-s},\ \widetilde{Q}_k\subset \Omega_{1/4}\},\ s=2,3,....
\end{equation}

\bigskip

\begin{lemma}\label{lf}
 If $q>n-1$, then
\begin{equation}\label{sum}
\sum\limits_{\widetilde{Q}_k\subset \Omega_{1/4}} d_k^q\leq C,
\end{equation}
where $C$ depends only on $n$, $q$ and $K$.
\end{lemma}

\begin{proof}
If $q\geq n$, \eqref{sum} is obvious. In the following, we consider the case when $q<n$.

For any  $Q_k\in \mathcal{F}^{s}$, there exists
$y_k\in (\partial\Omega)_1$ such that
$$\dist(Q_k,y_k)=\dist(Q_k,\partial\Omega_1)\leq 4d_k\leq2^{-s+2},$$
where Lemma \ref{l2.1} (iii) is used.
It follows from  $d_k\leq 2^{-s}$ for $Q_k\in \mathcal{F}^{s}$ that
$$\dist(x,(\partial\Omega)_1)\leq d_k+\dist(Q_k,y_k)
\leq 2^{-s}+2^{-s+2}\leq 2^{-s+3},\ \forall x\in Q_k\ \mathrm{and}\ Q_k\in \mathcal{F}^{s}$$
and then
\begin{equation}\label{ss1}
\mathcal{F}^{s}\subset\Omega_{1/4}\cap\{\dist(x,(\partial\Omega)_1)\leq 2^{-s+3}\}.
\end{equation}
By Lemma \ref{d}, we obtain
\begin{equation}\label{ss2}
|\mathcal{F}^{s}|\leq C{2^{-s}},
\end{equation}
where $C$ depends only on $n$ and $K$.

Observe that
$$\bigcup_{\widetilde{Q}_k\subset \Omega_{1/4}}Q_k
=\bigcup_{s=2}^{\infty}\bigcup_{Q_k\in \mathcal{F}^{s}} Q_k.
$$If $q>n-1$, we
derive from \eqref{fs} and \eqref{ss2} that
$$\begin{array}{rl}
\displaystyle\sum\limits_{\widetilde{Q}_k\subset \Omega_{1/4}} d_k^q
&\leq\displaystyle\sum\limits_{s=2}^{\infty}\left\{
\displaystyle\sum\limits_{Q_k\in \mathcal{F}^{s}}
\left(d_k^{q-n}\cdot d_k^{n}\right)\right\}
\leq \displaystyle\sum\limits_{s=2}^{\infty}\left\{2^{-s(q-n)}
\cdot\displaystyle\sum\limits_{Q_k\in \mathcal{F}^{s}}{d_k^{n}} \right\}\\
\\&\leq C \displaystyle\sum\limits_{s=2}^{\infty}2^{-s(q-n)}|\mathcal{F}^{s}|\leq
C \displaystyle\sum\limits_{s=2}^{\infty}2^{-s(q-n+1)}\leq C,
\end{array}$$
where $C$ depends only on $n$, $q$ and $K$.
\end{proof}

\bigskip

Fix $s_0\geq 2$ and a cube $Q_k\in \mathcal{F}^{s_0}.$ We classify the cubes $Q_l\in \mathcal{F}^s$ according to $\dist(Q_l,Q_k)$:
\begin{equation}\label{f}
\mathcal{F}^{s,j}_{Q_k}=\left\{
\begin{aligned}
&\bigcup_l\ \{Q_l\in \mathcal{F}^s,\ \dist(Q_l,Q_k)\leq 2^{-s_0+5}\}, &j=0,\\
&\bigcup_l\ \{Q_l\in \mathcal{F}^s,\ 2^{-s_0+j+4}<\dist(Q_l,Q_k)\leq 2^{-s_0+j+5}\},\ &j\geq 1.
\end{aligned}\right.
\end{equation}

First, we specify range of indexes $j$ and $s$ such that $\mathcal{F}^{s,j}_{Q_k}=\emptyset$. Since $\mathcal{F}^{s,j}_{Q_k}\subset \mathcal{F}^{s}\subset \Omega_{1/4}$, we see that $\mathcal{F}^{s,j}_{Q_k}=\emptyset$ for $j>s_0$. For any $Q_l\in \mathcal{F}^{s,j}_{Q_k}$, we have $Q_l\in \mathcal{F}^s$ and then
\begin{equation*}
\begin{aligned}
&2^{-s-1}<d_l\leq\dist(Q_l,\partial\Omega_1)\leq \dist(Q_l,Q_k)+\dist(Q_k,\partial\Omega_1)+\diam Q_k\\
&\mbox{}\hskip0.3cm= \dist(Q_l,Q_k)+\dist(Q_k,\partial\Omega_1)+d_k
\leq 2^{-s_0+j+5}+2^{-s_0+2}+2^{-s_0}
\leq 2^{-s_0+j+6},
\end{aligned}
\end{equation*}
where Lemma \ref{l2.1} (iii) is used.
Hence, $\mathcal{F}^{s,j}_{Q_k}=\emptyset$ for $s<s_0-j-6$.
In conclusion, for any fixed $s_0\geq 2$ and $Q_k\in \mathcal{F}^{s_0},$
\begin{equation}\label{o'}
\mathcal{F}^{s,j}_{Q_k}=\emptyset \ \ \mathrm{for} \ \ j> s_0\ \ \mathrm{or} \ \ s<s_0-j-6
\end{equation}
and by \eqref{3.11},
\begin{equation}\label{o}
\Omega_{1/{12}}\subset
\bigcup\limits_{\widetilde Q_l\subset \Omega_{1/4}}Q_l
=\bigcup_{j=0}^{s_0}\bigcup_{s=s_0-j-6}^{{\infty}}
\mathcal{F}^{s,j}_{Q_k}.
\end{equation}

\bigskip

\begin{lemma}\label{ff}
Fix $s_0\geq 2$ and a cube $Q_k\in \mathcal{F}^{s_0}.$ There exists a constant $C$ depending only on $n$ and $K$ such that
\begin{equation}\label{m}
|\mathcal{F}^{s,j}_{Q_k}|\leq C2^{(-s_0+j)(n-1)-s}.
\end{equation}
\end{lemma}
\begin{proof}
By \eqref{o'}, we only consider the case when $0\leq j\leq s_0$ and $s\geq s_0-j-6$. Since $Q_k\in \mathcal{F}^{s_0}\subset \Omega_{1/4}$, there exists $y_k\in (\partial\Omega)_{1/2}$ such that
$$\dist(Q_k,y_k)=\dist(Q_k,\partial\Omega_1)\leq 4d_k\leq2^{-s_0+2},$$
where Lemma \ref{l2.1} (iii) is used.
It follows that for any $Q_l\in \mathcal{F}^{s,j}_{Q_k}$,
$$\begin{array}{l}\dist(Q_l,y_k)\leq \dist(Q_l,Q_k)+\dist(Q_k,y_k)+\diam Q_k\\[10pt]
\mbox{}\hskip.3cm=
\dist(Q_l,Q_k)+\dist(Q_k,y_k)+d_k
\leq
2^{-s_0+j+5}+2^{-s_0+2}+2^{-s_0}\leq 2^{-s_0+j+6}
\end{array}
$$
and by $d_l\leq 2^{-s}\leq 2^{-s_0+j+6}$,
$$\dist(x,y_k)\leq d_l+\dist(Q_l,y_k)\leq 2^{-s_0+j+6}+2^{-s_0+j+6}\leq 2^{-s_0+j+7},\ \forall x\in Q_l.$$
Combing the above estimate and \eqref{ss1}, we obtain
$$\mathcal{F}^{s,j}_{Q_k}\subset \Omega_{2^{-s_0+j+7}}(y_k)\cap\{\dist(x,(\partial\Omega)_1)\leq 2^{-s+3}\}.$$
\eqref{m} then follows by Lemma \ref{d}.
\end{proof}

\bigskip

\begin{lemma}\label{g}
Fix a cube $Q_l\in \mathcal{F}^{s}$ and change $Q_k\in \mathcal{F}^{s_0}$. Then
there exist at most $C2^{j(n-1)}$ cubes $Q_k\in\mathcal{F}^{s_0}$ such that $Q_l\in  \mathcal{F}^{s,j}_{Q_k}$, where $C$ depends only on $n$ and $K$.
\end{lemma}
\begin{proof}
By \eqref{o'}, we only consider the case when  $0\leq j\leq s_0$ and $s\geq s_0-j-6$. If $Q_{l}\in \mathcal{F}^{s,j}_{Q_k}$, then $\dist(Q_l,Q_k)\leq 2^{-s_0+j+5}$. It then follows from $d_k\leq 2^{-s_0}$
and $d_l\leq 2^{-s}\leq 2^{-s_0+j+6}$ that
$$\dist(x,x_l)\leq d_k+\dist(Q_l,Q_k)+d_l\leq  2^{-s_0+j+7},\ \forall x\in Q_k,$$
where $x_l$ is the center of $Q_l$.
We deduce from \eqref{ss1} and the above estimate that
$$\displaystyle\bigcup_{k}\ \{Q_k: Q_k\in \mathcal{F}^{s_0}\ s.t.\ Q_{l}\in \mathcal{F}^{s,j}_{Q_k}\}\subset \{|x'-x_{l}'|\leq 2^{-s_0+j+7},\ \dist(x,(\partial\Omega)_1)\leq 2^{-s_0+3}\}$$
and then by Lemma \ref{d},
$$\bigg|\displaystyle\bigcup_{k}\ \{Q_k: Q_k\in \mathcal{F}^{s_0}\ s.t.\ Q_{l}\in \mathcal{F}^{s,j}_{Q_k}\}\bigg|\leq C2^{(-s_0+j)(n-1)-s_0}=C2^{j(n-1)-ns_0}, $$
where $C$ is a constant depending only on $n$ and $K$.
Since for any $Q_k\in\mathcal{F}^{s_0}$, $|Q_k|=C_nd_k^n\geq C_n2^{-ns_0-n}$ for some constant $C_n$ depending only on $n$,  we have
$\bigcup_{k}\ \{Q_k: Q_k\in \mathcal{F}^{s_0}\ s.t.\ Q_{l}\in \mathcal{F}^{s,j}_{Q_k}\}$ contains at most $\frac{2^nC}{C_n }2^{j(n-1)}$ cubes in $\mathcal{F}^{s_0}$.
The lemma is thus proved.
\end{proof}

\section{Preliminary results}

We start with the following $L^p$ estimates.

\bigskip

\begin{lemma}\label{lb}
Let $\Omega$ be a Lipschitz domain in $\mathbb{R}^n$ and $u$ satisfy
\begin{equation*}
\left\{
\begin{array}{lr}
\Delta u=f\ \ \ \mathrm{in}\ \ \ \Omega,\\ \ \ u=0\ \ \ \mathrm{on}\ \ \partial\Omega,
\end{array}
\right.
\end{equation*}
where $$f\in L^p(\Omega) \ \ \mathrm{for}\ \ 1<p<{\infty} \ \ \mathrm{and}\ \ \supp f\subset D\subset\Omega.$$
Then  we have for any measurable set $E\subset \Omega$,
\begin{equation}\label{xxj}
||u||_{L^p(E)}\leq C|E|^\frac{2}{np}|D|^\frac{2}{np'}||f||_{L^p(D)},
\end{equation}
where $C$ depends only on $n$ and $p$.
\end{lemma}
\begin{proof}
Let $G=G(x,y)$ be the (Dirichlet) Green's function of the operator $-\Delta$ on the domain $\Omega$. By Green's representation formula, we have
$$u(x)=\int_{\Omega}G(x,y)f(y)dy,\ \ \forall x\in \Omega.$$
Let $\Gamma=\Gamma(x,y)$ be the normalized fundamental solution of Laplace's equation. For $n\geq 3$, by comparison principle,
$$0\leq G(x,y)\leq \Gamma(x,y)= C_n|x-y|^{2-n},\ \ \forall x,y\in\Omega,$$
where $C_n$ depends only on $n$.
Hence,
$$|u(x)|\leq C_n\int_{D}|x-y|^{2-n}|f(y)|dy,\ \ \forall x\in \Omega.$$

By H\"{o}lder's inequality,
\begin{equation*}
\begin{aligned}
\displaystyle\int_{D}|x-y|^{2-n}|f(y)|dy
&=\displaystyle\int_{D}|x-y|^{\frac{2-n}{p'}}|x-y|^{\frac{2-n}{p}}|f(y)|dy\\
&\leq \left(\int_{D}|x-y|^{2-n}dy\right)^{\frac{1}{p'}}
\left(\int_{D}|x-y|^{2-n}|f(y)|^pdy\right)^{\frac{1}{p}}.
\end{aligned}
\end{equation*}
Choose $R>0$ such that $|D|=|B_R(x)|$ and we deduce that
\begin{equation}\label{de}
\begin{aligned}
\int_{D}|x-y|^{2-n}dy\leq \int_{B_R(x)}|x-y|^{2-n}dy
\leq C_nR^{2}\leq C_n|D|^{\frac{2}{n}},\ \ \forall x\in \mathbb{R}^n.
\end{aligned}
\end{equation}
Therefore,
$$|u(x)|\leq C|D|^{\frac{2}{np'}}\left(\int_{D}|x-y|^{2-n}|f(y)|^pdy\right)^{\frac{1}{p}},$$
where $C$ depends only on $n$ and $p$.
It follows that
\begin{equation*}
\begin{aligned}
\displaystyle\int_{E}|u(x)|^{p}dx
&\leq C|D|^{\frac{2p}{np'}}\int_{E}\int_{D}|x-y|^{2-n}|f(y)|^pdydx\\
&\leq C|D|^{\frac{2p}{np'}}
\int_{D}|f(y)|^pdy
\sup\limits_{y\in D}\displaystyle\int_{E}|x-y|^{2-n}dx.
\end{aligned}
\end{equation*}
Similar to the derivation of \eqref{de}, we have
$$\sup\limits_{y\in D}\displaystyle\int_{E}|x-y|^{2-n}dx\leq C|E|^{\frac{2}{n}}$$
and then
$$\displaystyle\int_{E}|u(x)|^{p}dx\leq
C |E|^{\frac{2}{n}} |D|^{\frac{2p}{np'}}
\int_{D}|f(y)|^pdy,$$
which implies \eqref{xxj} for $n\geq 3$.

The proof for $n=2$ is similar and we omit it here.
\end{proof}

\bigskip

The next lemma concerns pointwise boundary  $C^{1,\alpha}$ estimates and we refer to Theorem 1.6 in \cite{LZ} for its proof.

\bigskip

\begin{lemma}\label{l3.2}
Assume $0\in \partial\Omega$ and there exists $\varphi\in C^{1,\alpha}(B_1')$ such that
$$\Omega_{1}=B_{1}\cap\{ x^n>\varphi(x')\}\ \ \mathrm{and}\ \ (\partial\Omega)_{1}=B_{1}\cap\{ x^n=\varphi(x')\}$$ for $0<\alpha<1$. Let $u$ satisfy
\begin{equation*}
\left\{
\begin{array}{lr} \Delta u=f\ \ \ \mathrm{in}\ \ \ \ \Omega_1,\\ \ \ u=g\ \ \ \mathrm{on}\ \ (\partial\Omega)_1,
\end{array}
\right.
\end{equation*}
where $g\in C^{1,\alpha}(0)$ and $f\in L^n(\Omega_1)$ such that
\begin{equation*}
||f||_{L^n(\Omega_r)}\leq K_f r^{\alpha}, \ \forall 0<r<1
\end{equation*}
for some constant $K_f$.
Then $u\in C^{1,\alpha}(0)$, i.e., there exists an affine function $l$ such that
$$|u(x)-l(x)|\leq C|x|^{1+\alpha}
(||u||_{L^{\infty}(\Omega_1)}+K_f+||g||_{C^{1,\alpha}(0)}),\ \forall x\in \Omega_{r_0}$$
and
$$|Dl|\leq C(||u||_{L^{\infty}(\Omega_1)}+K_f+||g||_{C^{1,\alpha}(0)}),$$
where $C$ and $r_0$ depends on $n$, $\alpha$ and  $||\varphi||_{C^{1,\alpha}(B_1')}$.
\end{lemma}

\bigskip
\begin{re}\label{r3.3}
If $f\in L^p(\Omega_1)$ for $p>n$, then by H$\ddot{o}$lder's inequality, we have
$$||f||_{L^n(\Omega_r)}\leq C_n r^{1-n/p}, \ \forall 0<r<1,$$
where $C_n$ depends only on $n$. From Lemma \ref{l3.2}, we obtain pointwise boundary $C^{1,\min\{\alpha,1-n/p\}}$ regularity, which is optimal by Sobolev embedding theorem.
\end{re}

\bigskip

\begin{cor}\label{c3.0}
Let $u$ satisfy
\begin{equation*}
\left\{
\begin{array}{lr} \Delta u=0\ \ \ \mathrm{in}\ \ \ \ \Omega_r,\\ \ \ u=0\ \ \ \mathrm{on}\ \ (\partial\Omega)_r
\end{array}
\right.
\end{equation*}
with $0<r\leq 1$.
Then  $u$ is $C^{1,\alpha}$ at $x_0$ for any $x_0\in (\partial\Omega)_{r/2}$, i.e., there exists an affine function $l_{x_0}$ such that for any $1<p<\infty,$
\begin{equation}\label{lzz}
|u(x)-l_{x_0}(x)|\leq  C r^{-(1+\alpha+n/p)} |x-x_0|^{1+\alpha}||u||_{L^p(\Omega_r)},\ \forall x\in \Omega_{3r/4}
\end{equation}
and
\begin{equation}\label{gl}
|Dl_{x_0}|\leq C r^{-1-n/p}||u||_{L^p(\Omega_r)},
\end{equation}
where $C$ depends on $n,\alpha,p$ and $||\varphi||_{C^{1,\alpha}(B_r')}$.
\end{cor}
\begin{proof}
We may assume that $r=1$. By boundary local maximum principle (see Theorem 9.26 in \cite{GT}), we have
$||u||_{L^\infty(\Omega_{3/4})}\leq C||u||_{L^{p}(\Omega_1)}$, where $C$ depends only on $n$ and $p$. Then \eqref{lzz} and \eqref{gl} follow from Lemma \ref{l3.2} and standard scaling arguments.
\end{proof}

\bigskip

We end this section by the following local boundary $W^{2,p}$ estimates for harmonic functions on $C^{1,\alpha}$ domains.

\bigskip

\begin{thm}\label{t3.1}
Let $1<p<{\infty}$ and $1-1/p<\alpha\leq1$. Assume that $0\in \partial\Omega$ and there exists $\varphi\in C^{1,\alpha}(B_1')$ such that
$$\Omega_{1}=B_{1}\cap\{ x^n>\varphi(x')\}\ \ \mathrm{and}\ \ (\partial\Omega)_{1}=B_{1}\cap\{ x^n=\varphi(x')\}.$$ If $u\in W^{2,p}(\Omega_1)$ satisfies
\begin{equation*}
\left\{
\begin{array}{lr} \Delta u=0\ \ \ \mathrm{in}\ \ \ \ \Omega_1,\\ \ \ u=0\ \ \ \mathrm{on}\ \ (\partial\Omega)_1,
\end{array}
\right.
\end{equation*}
then we have
\begin{equation}\label{h}
||D^2 u||_{ L^{p}(\Omega_{1/12})}\leq C||u||_{L^{p}(\Omega_1)},
\end{equation}
where $C$ depends on $n,\alpha,p$ and $||\varphi||_{C^{1,\alpha}(B_1')}$.
\end{thm}
\begin{proof}
Let $\{Q_k\}_{k=1}^{\infty}$ be Whitney decomposition of $\Omega_{1}$ and ${\widetilde Q_k}=\frac65 Q_k$.
For any $\widetilde Q_k\subset \Omega_{1/4}$, we let $y_k\in (\partial\Omega)_{1/2}$ and $\tilde x_k\in \partial\widetilde Q_k$ such that
$$|\tilde x_k-y_k|=\dist(\widetilde Q_k, \partial\Omega_1)<\dist(Q_k, \partial\Omega_1)\leq 4d_k,$$
where Lemma \ref{l2.1} (iii) is used in the last inequality.
Consequently, we see that
$$ |x-y_k|\leq|x-\tilde x_k|+|\tilde x_k-y_k|\leq 6d_k,\ \forall x\in\widetilde Q_k.$$
It then follows from  Corollary \ref{c3.0} that $u\in C^{1,\alpha}(y_k)$ and  there exists an affine function $l_{y_k}$(written by $l$ for simplicity in the following) such that
\begin{equation*}
|(u-l)(x)|\leq  C |x-y_k|^{1+\alpha}||u||_{L^p(\Omega_{1})}
\leq C d_k^{1+\alpha}||u||_{L^p(\Omega_{1})},
\ \forall x\in \widetilde Q_k,
\end{equation*}
where $C$ depends on $n,\alpha, p$ and $||\varphi||_{C^{1,\alpha}(B_1')}$.

Since
$u-l$ satisfies
$\Delta(u-l)=0$ in $\widetilde Q_k,$
we deduce from interior $W^{2,p}$ estimate and the above pointwise $C^{1,\alpha}$ estimate that
\begin{equation}\label{3.12}
||D^2 (u-l)||_{L^p(Q_k)}
\leq C d_k^{n/p-2}||u-l||_{L^{\infty}(\widetilde Q_k)}
\leq C d_k^{n/p+\alpha-1}||u||_{L^p(\Omega_{1})},
\end{equation}
where  $C$ depends on $n,\alpha,p$ and $||\varphi||_{C^{1,\alpha}(B_1')}$.

In conclusion, since
 $n-(1-\alpha)p>n-1$ as $1/(1-p)<\alpha\leq 1,$
we infer from  Lemma \ref{l2.3}, \eqref{3.12}  and  Lemma \ref{lf} that
\begin{equation*}
\begin{aligned}
||D^2 u||_{L^p(\Omega_{1/{12}})}^{p}
\leq
&\sum\limits_{\widetilde Q_k\subset \Omega_{1/4}}||D^2 u||_{L^p(Q_k)}^{p}
\\ \leq C ||u||_{L^p(\Omega_{1})}^p
&\displaystyle\sum\limits_{\widetilde Q_k\subset \Omega_{1/4}}
d_k^{n-(1-\alpha)p}
\leq C||u||_{L^p(\Omega_{1})}^p,
\end{aligned}
\end{equation*}
where $C$ depends on $n,\alpha,p$ and $||\varphi||_{C^{1,\alpha}(B_1')}$.
\end{proof}

\bigskip

\begin{re}\label{hi}
Theorem 3.4 follows from
interior $W^{2,p}$ estimate, boundary $C^{1,\alpha}$ estimate (Corollary \ref{c3.0}) and
Whitney decomposition (Lemma 2.1, 2.2 and \ref{lf}).
If we apply this argument to non-homogeneous equation $\Delta u=f$, then $p>n$ is needed and we can only arrive at
$$||D^2 u||_{ L^{p_0}(\Omega_{1/12})}\leq C(||u||_{L^{p}(\Omega_1)}+||f||_{L^p(\Omega_1)})\ \ \mathrm{with}\ \ 1\leq p_0<\min\{1/(1-\alpha),p/n\}$$
since only $u\in C^{1,\min\{\alpha,1-\frac{n}{p}\}}$ can be obtained as $p>n$ (cf. Remark \ref{r3.3}).
To improve the above estimate and remove the restriction $p>n$, we need to decompose $u$ according to Whitney decomposition (see more details in Section 4).

\end{re}

\section{$W^{2,p}$ estimate for Poisson's equation on $C^{1,\alpha}$ domain}

By considering $u-g$ and using the technique of perturbation from the constant coefficient case, Theorem \ref{t1.1} follows easily from the following  $W^{2,p}$ estimates for Poisson's equation.

\bigskip

\begin{thm}\label{t4.1}
Let $1<p<{\infty}$ and $1-1/p<\alpha\leq1$. Assume $0\in \partial\Omega$ and there exists $\varphi\in C^{1,\alpha}(B_1')$ such that
$$\Omega_{1}=B_{1}\cap\{ x^n>\varphi(x')\}\ \ \mathrm{and}\ \ (\partial\Omega)_{1}=B_{1}\cap\{ x^n=\varphi(x')\}.$$
If $u\in W^{2,p}(\Omega_1)$ and $f\in L^{p}(\Omega_1)$ such that
\begin{equation}\label{4.1}
\left\{
\begin{array}{lr} \Delta u=f\ \ \ \mathrm{in}\ \ \ \ \Omega_1,\\ \ \ u=0\ \ \ \mathrm{on}\ \ \ (\partial\Omega)_1,
\end{array}
\right.
\end{equation}
 then we have
\begin{equation}\label{u}
||D^2 u||_{ L^{p}(\Omega_{1/24})}\leq C(||u||_{L^{p}(\Omega_1)}+||f||_{L^{p}(\Omega_1)}),
\end{equation}
where $C$ depends on $n,\alpha,p$ and $||\varphi||_{C^{1,\alpha}(B_1')}$.
\end{thm}

\begin{proof}
Let $\{Q_l\}_{l=1}^{\infty}$ be Whitney decomposition of $\Omega_1$ and $\widetilde{Q}_l=\frac65{Q_l}$. We separate $u$ to be $$u=v+w$$ such that
\begin{equation*}
\left\{{
\begin{array}{l} {\Delta v=f\chi_{\cup_{\widetilde Q_l\not\subset\Omega_{1/4}}}\ \mathrm{in}\ \ \ \Omega_1,}\\
{ \ \ v=u\ \ \ \ \ \ \ \ \ \ \ \ \ \ \mathrm{on} \ \ \partial \Omega_1}
\end{array}}
\right.
\ \ \mathrm{and}\ \ \
\left\{ {
\begin{array}{l}{\Delta w=f\chi_{\cup_{\widetilde Q_l\subset\Omega_{1/4}}}\ \ \mathrm{in}\ \ \ \ \Omega_1},\\
{ \ \ w=0\ \ \ \ \ \ \ \ \ \ \ \ \  \ \ \mathrm{on}\ \ \ \partial\Omega_1.}
\end{array}}
\right.
\end{equation*}

Since, by Lemma \ref{l2.3},
$\Omega_{1/{12}}\subset\bigcup_{\widetilde Q_k\subset \Omega_{1/{4}}} Q_k$,
we have $v$ is harmonic in $\Omega_{1/{12}}$. And then by Theorem \ref{t3.1},
\begin{equation}\label{v}
||D^2 v||_{ L^{p}(\Omega_{1/{24}})}\leq C||v||_{L^{p}(\Omega_{1/{12}})}\leq C(||u||_{L^{p}(\Omega_1)}+||f||_{L^{p}(\Omega_1)}).
\end{equation}

Our sequent work is devoted to prove the following estimate:
\begin{equation}\label{w}
||D^2 w||_{ L^{p}(\Omega_{1/{12}})}\leq C||f||_{L^{p}(\Omega_1)},
\end{equation}
where $C$ depends on $n,\alpha,p$ and $||\varphi||_{C^{1,\alpha}(B_1')}$.

For this purpose, we decompose $w$ according to Whitney decomposition as follows. Set
$$\mathcal{F}^{s}=\bigcup_{k}\ \{Q_k:2^{-s-1}< d_k\leq 2^{-s},\ \widetilde{Q}_k\subset \Omega_{1/4}\},\ s=2,3,...$$
as in  \eqref{fs}.
Fix $s_0\geq 2$ and  $Q_k\in\mathcal{F}^{s_0}$. Let
\begin{equation*}
\mathcal{F}^{s,j}_{Q_k}=\left\{
\begin{aligned}
&\bigcup_l\ \{Q_l\in \mathcal{F}^s,\ \dist(Q_l,Q_k)\leq 2^{-s_0+5}\}, &j=0,\\
&\bigcup_l\ \{Q_l\in \mathcal{F}^s,\ 2^{-s_0+j+4}<\dist(Q_l,Q_k)\leq 2^{-s_0+j+5}\},\ &j\geq 1
\end{aligned}\right.
\end{equation*}
as in \eqref{f} and $w^{s,j}_{k}$ satisfy
\begin{equation}\label{equ}
\left\{
\begin{array}{lr} \Delta w^{s,j}_{k}=f\chi_{{F}^{s,j}_{Q_k}}\ \ \ \mathrm{in}\ \ \ \ \Omega_1,\\ \ \ \ w^{s,j}_{k}=0\ \ \ \ \ \  \ \ \ \mathrm{on}\ \ \ \partial\Omega_1.
\end{array}
\right.
\end{equation}
Recall \eqref{o}, that is,
$$\Omega_{1/{12}}\subset\bigcup\limits_{\widetilde Q_l\subset \Omega_{1/4}}Q_l
=\bigcup_{j=0}^{s_0}\bigcup_{s=s_0-j-6}^{{\infty}}
\mathcal{F}^{s,j}_{Q_k}.$$
It then follows that
\begin{equation}\label{wsum}
w=\sum\limits_{j=0}^{s_0}
\sum\limits_{s=s_0-j-6}^{{\infty}}w^{s,j}_{k}\ \ \mathrm{in}\ \ \Omega_1.
\end{equation}

We need the following estimate of $w^{s,j}_{k}$.

\bigskip

\begin{lemma}\label{l4.2}
Fix $s_0\geq 2$ and  $Q_k\in\mathcal{F}^{s_0}$.
Let
\begin{equation}\label{mbeta}
m=s-s_0\ \ \mathrm{and} \ \ \beta=\alpha+\frac{n}{p}+\frac{2}{np'}-1.
\end{equation}
Then
\begin{equation}\label{wsjk}||D^2 w^{s,j}_{k}||_{L^{p}(Q_k)}
\leq C 2^{-j\beta-\frac{2m}{np'}}
||f||_{L^p(\mathcal{F}^{s,j}_{Q_k})},
\end{equation}
where $C$ depends on $n,\alpha,p$ and $||\varphi||_{C^{1,\alpha}(B_1')}$.
\end{lemma}

\begin{proof}
We divide the proof of  \eqref{wsjk} into two cases: $j=0$ and $j\geq1$.

(i) As $j=0$, by \eqref{equ} and Lemma \ref{lb}, we have
$$|| w^{s,0}_{k}||_{L^p(\widetilde Q_k)}
\leq C|\widetilde Q_k|^{\frac{2}{np}}|\mathcal{F}^{s,0}_{Q_k}|^{\frac{2}{np'}}
||f||_{L^p(\mathcal{ F}^{s,0}_{Q_k})}.$$
Since $Q_k\in \mathcal{F}^{s_0}$, $|\widetilde Q_k|\leq 2^{-s_0n}$. By Lemma \ref{ff},
$|\mathcal{F}^{s,0}_{Q_k}|\leq C2^{-s_0(n-1)-s}$. Thus,
\begin{equation*}
\begin{aligned}
|| w^{s,0}_{k}||_{L^p(\widetilde Q_k)}
\leq C2^{-2s_0-\frac{2m}{np'}}||f||_{L^p(\mathcal{F}^{s,0}_{Q_k})}.
\end{aligned}
\end{equation*}
In view of $\Delta w^{s,0}_{k}=f\chi_{\mathcal{F}^{s,0}_{Q_k}}$ in $\widetilde Q_k$, we deduce
from interior $W^{2,p}$ estimate that
\begin{equation}\label{00}
\begin{aligned}
||D^2 w^{s,0}_{k}||_{L^{p}(Q_k)}
&\leq C\left(d_k^{-2}|| w^{s,0}_{k}||_{L^{p}(\widetilde Q_k)}
+ ||f\chi_{\mathcal{F}^{s,0}_{Q_k}}||_{L^{p}(\widetilde Q_k)}\right)\\
\\&\leq C\left(2^{-\frac{2m}{np'}}||f||_{L^p(\mathcal{F}^{s,0}_{Q_k})}+
||f||_{L^p(\mathcal{F}^{s,0}_{Q_k}\cap \widetilde Q_k)}\right).
\end{aligned}
\end{equation}

If $m<-6$, we infer from \eqref{o'} that $\mathcal{F}^{s,j}_{Q_k}=\emptyset$. Hence we only need consider $m\geq-6$.
As $-6\leq m<5$, \eqref{wsjk} follows easily from \eqref{00}.
As $m\geq 5$, we claim $\mathcal{F}^{s,0}_{Q_k}\cap\widetilde Q_k=\emptyset$. Actually,
for any $Q_l\in \mathcal{F}^{s,0}_{Q_k}$, we have
$$\dist(x,\partial\Omega_1)\leq d_l+\dist(Q_l,\partial\Omega_1)\leq 5d_l\leq 2^{-s+3}=2^{-m-s_0+3}\leq2^{-s_0-2},\ \forall x\in Q_l.$$
However, since $Q_k\in \mathcal{F}^{s_0}$,
$$\dist(\widetilde Q_k,\partial\Omega_1)
\geq\dist(Q_k,\partial\Omega_1)-d_k/5\geq4d_k/5> 2^{-s_0-2}.$$
Hence $\mathcal{F}^{s,0}_{Q_k}\cap\widetilde Q_k=\emptyset$ and then we derive \eqref{wsjk} from \eqref{00}.~\\

(ii) As $j\geq1$, let $y_k\in (\partial\Omega)_{1/2}$ and  $z_k\in Q_k$ such that
$|z_k-y_k|=\dist(Q_k,\partial\Omega)$ and we claim the following relation which will be used several times:
\begin{equation}\label{11.22}
\widetilde Q_k\subset\Omega_{2^{-s_0+3}}(y_k) \subset \Omega_{2^{-s_0+j+3}}(y_k)\subset\{\dist(x,Q_k)\leq 2^{-s_0+j+4}\},
\end{equation}
where $\Omega_{2^{-s_0+j+3}}(y_k)=\Omega\cap B_{2^{-s_0+j+3}}(y_k).$
Indeed, we infer from $Q_k\in \mathcal{F}^{s_0}$ that $d_k\leq 2^{-s_0}$ and then by Lemma \ref{l2.1} (iii),
$$|z_k-y_k|=\dist(Q_k,\partial\Omega)\leq4d_k\leq 2^{-s_0+2}.$$
Since $\diam \widetilde  Q_k=\frac65d_k\leq 2^{-s_0+1}$, we have
$$|x-y_k|\leq |x-z_k|+|z_k-y_k|\leq \diam \widetilde  Q_k+2^{-s_0+2}\leq 2^{-s_0+3},\ \forall x\in \widetilde Q_k,$$
which implies that $\widetilde Q_k\subset\Omega_{2^{-s_0+3}}(y_k)$.
Since
$$\dist(x,Q_k)\leq |x-z_k|\leq |x-y_k|+|y_k-z_k|\leq 2^{-s_0+j+3}+2^{-s_0+2} \leq2^{-s_0+j+4}$$ for any $x\in\Omega_{2^{-s_0+j+3}}(y_k)$, we have
$\Omega_{2^{-s_0+j+3}}(y_k)\subset\{\dist(x,Q_k)\leq 2^{-s_0+j+4}\}$ and then \eqref{11.22} holds.

By the definition of $\mathcal{F}^{s,j}_{Q_k}$, we have $\dist(Q_l,Q_k)>2^{-s_0+j+4}$ as
$Q_l\in \mathcal{F}^{s,j}_{Q_k}$ and then  $$\mathcal{F}^{s,j}_{Q_k}\subset \{\dist(x,Q_k)>2^{-s_0+j+4}\}.$$
Combining it with \eqref{11.22}, we get $\mathcal{F}^{s,j}_{Q_k}\cap \Omega_{2^{-s_0+j+3}}(y_k)=\emptyset$ and then by \eqref{equ},
\begin{equation}\label{j2}
\Delta{w}^{s,j}_{k}=0\ \ \mathrm{in}\ \ \Omega_{2^{-s_0+j+3}}(y_k).
\end{equation}
From Corollary \ref{c3.0}, it follows that
there exists an affine function $l$ such that for any $x\in \Omega_{2^{-s_0+j+2}}(y_k)$,
$$
|({w}^{s,j}_{k}-l)(x)|
\leq C 2^{-(j-s_0)(1+\alpha+n/p)}|x-y_k|^{1+\alpha}
||{w}^{s,j}_{k}||_{L^p(\Omega_{2^{-s_0+j+3}}(y_k))}.
$$

From \eqref{equ} and Lemma \ref{lb}, we deduce that
\begin{equation*}
\begin{aligned}
||\displaystyle {w}^{s,j}_{k}||_{L^p(\Omega_{2^{-s_0+j+3}}(y_k))}
&\leq C|B_{2^{-s_0+j+3}}|^{\frac{2}{np}}|\mathcal{F}^{s,j}_{Q_k}|^{\frac{2}{np'}}
||f||_{L^p(\mathcal{F}^{s,j}_{Q_k})}\\
&\leq C2^{-2(s_0-j-\frac{s_0-s-j}{np'})}||f||_{L^p(\mathcal{F}^{s,j}_{Q_k})},
\end{aligned}
\end{equation*}
where $|\mathcal{F}^{s,j}_{Q_k}|\leq C2^{(-s_0+j)(n-1)-s}$ is used that is given by Lemma \ref{ff}.

By \eqref{11.22}, $\widetilde Q_k\subset\Omega_{2^{-s_0+3}}(y_k)$ and then for any $x\in\widetilde Q_k$,
$$|x-y_k|\leq 2^{-s_0+3}.$$

Combining above estimates, we obtain
\begin{equation*}
\begin{aligned}
||{w}^{s,j}_{k}-l||_{L^\infty(\widetilde Q_k)}
&\leq C 2^{-(j-s_0)(1+\alpha+n/p)-s_0(1+\alpha)}
||{w}^{s,j}_{k}||_{L^p(\Omega_{2^{-s_0+j+3}}(y_k))}
\\
&\leq C 2^{
-s_0(2-\frac np)-j\beta-\frac{2m}{np'}}
||f||_{L^p(\mathcal{F}^{s,j}_{Q_k})},
\end{aligned}
\end{equation*}
where $\beta=\alpha+\frac{n}{p}+\frac{2}{np'}-1$ is defined by \eqref{mbeta}.

In view of \eqref{11.22} and \eqref{j2},
$$\Delta({w}^{s,j}_{k}-l)=0\ \ \mathrm{in}\ \ \widetilde Q_k.$$
Using interior $W^{2,p}$ estimate,
\begin{equation*}
\begin{array}{ccc}
&||D^2 w^{s,j}_{k}||_{L^{p}(Q_k)}
\leq C d_k^{n/p-2}||{w}^{s,j}_{k}-l||_{L^{\infty}(\widetilde Q_k)}
\leq C 2^{-j\beta
-\frac{2m}{np'}}
||f||_{L^p(\mathcal{F}^{s,j}_{Q_k})},
\end{array}
\end{equation*}
where $C$ depends on $n,\alpha,p$ and $||\varphi||_{C^{1,\alpha}(B_1')}$.
\end{proof}

\bigskip

Now we continue the proof of Theorem \ref{t4.1}.

By Lemma \ref{l2.3}, we deduce that
$$\Omega_{1/{12}}\subset\bigcup_{\widetilde Q_k\subset \Omega_{1/{4}}} Q_k
= \bigcup_{s_0=2}^{\infty}
 \bigcup_{ Q_k\in\mathcal{F}^{s_0}}Q_k$$
and then
$$||D^2 w||_{L^p(\Omega_{1/{12}})}^{p}
\leq
\sum\limits_{\widetilde Q_k\subset \Omega_{1/{4}}}||D^2 w||_{L^p(Q_k)}^{p}
= \sum\limits_{s_0=2}^{\infty}
\sum\limits_{Q_k\in\mathcal{F}^{s_0}}||D^2 w||_{L^p(Q_k)}^{p}.$$
From \eqref{wsum} and  Minkowski's inequality, it follows that
$$||D^2 w||_{L^p(Q_k)}\leq \sum\limits_{j=0}^{s_0}
\sum\limits_{s=s_0-j-6}^{{\infty}}
||D^2 w^{s,j}_{k}||_{L^p(Q_k)}$$
and then
\begin{equation}\label{h}
\begin{aligned}||D^2 w||_{L^p(\Omega_{1/{12}})}^{p}
\leq \sum\limits_{s_0=2}^{\infty}
\sum\limits_{Q_k\in\mathcal{F}^{s_0}}
\left(\sum\limits_{j=0}^{s_0}
\sum\limits_{s=s_0-j-6}^{{\infty}}
||D^2 w^{s,j}_{k}||_{L^p(Q_k)}\right)^p.
\end{aligned}
\end{equation}

Let $\tau>0$ (depending on $n,\alpha$ and $p$) to be determined later and by H$\rm{\ddot o}$lder's inequality,
\begin{equation*}
\begin{aligned}
&\left(\sum\limits_{j=0}^{s_0}
\sum\limits_{s=s_0-j-6}^{{\infty}}
||D^2 w^{s,j}_{k}||_{L^p(Q_k)}\right)^p
\leq C
\sum\limits_{j=0}^{s_0}\left\{2^{j\tau p}\left(\sum\limits_{s=s_0-j-6}^{{\infty}}
||D^2 w^{s,j}_{k}||_{L^p(Q_k)}\right)^p\right\}\\
&\leq C \sum\limits_{j=0}^{s_0}\left\{2^{j\tau p}\left(\sum\limits_{s=s_0-j-6}^{{s_0}}
||D^2 w^{s,j}_{k}||_{L^p(Q_k)}\right)^p+
2^{j\tau p}\left(\sum\limits_{s=s_0+1}^{{\infty}}
||D^2 w^{s,j}_{k}||_{L^p(Q_k)}\right)^p\right\}.
\end{aligned}
\end{equation*}
Using H$\rm{\ddot o}$lder's inequality again,
$$\left(\sum\limits_{s=s_0-j-6}^{{s_0}}
||D^2 w^{s,j}_{k}||_{L^p(Q_k)}\right)^p\leq
C
\sum\limits_{s=s_0-j-6}^{{s_0}}
2^{(s_0-s)\tau p}||D^2 w^{s,j}_{k}||^p_{L^p(Q_k)}$$
and
$$\left(\sum\limits_{s=s_0+1}^{\infty}
||D^2 w^{s,j}_{k}||_{L^p(Q_k)}\right)^p\leq
C
\sum\limits_{s=s_0+1}^{\infty}
2^{(s-s_0)\tau p}||D^2 w^{s,j}_{k}||^p_{L^p(Q_k)}.$$
Recall $m=s-s_0$ given by \eqref{mbeta}. We derive from the above estimates that
$$\left(\sum\limits_{j=0}^{s_0}
\sum\limits_{s=s_0-j-6}^{{\infty}}
||D^2 w^{s,j}_{k}||_{L^p(Q_k)}\right)^p\leq C
\sum\limits_{j=0}^{s_0}
\sum\limits_{s=s_0-j-6}^{\infty}2^{(j+|m|)\tau p}
||D^2 w^{s,j}_{k}||_{L^p(Q_k)}^p.$$
Substitute it into \eqref{h} and consequently,
$$||D^2 w||_{L^p(\Omega_{1/{12}})}^{p}
\leq C \sum\limits_{s_0=2}^{\infty}
\sum\limits_{Q_k\in\mathcal{F}^{s_0}}\sum\limits_{j=0}^{s_0}
\sum\limits_{s=s_0-j-6}^{\infty}2^{(j+|m|)\tau p}
||D^2 w^{s,j}_{k}||_{L^p(Q_k)}^p.$$

By exchanging summation order,
\begin{equation*}
\begin{aligned}
||D^2 w||_{L^p(\Omega_{1/{12}})}^{p}
\leq C \sum\limits_{s=-4}^{\infty}\sum\limits_{j=0}^{\infty}
\sum\limits_{s_0=2}^{s+j+6}
\sum\limits_{Q_k\in\mathcal{F}^{s_0}}
2^{(j+|m|)\tau p}
||D^2 w^{s,j}_{k}||_{L^p(Q_k)}^p.
\end{aligned}
\end{equation*}
From Lemma \ref{l4.2}, it follows that
$$||D^2 w||_{L^p(\Omega_{1/{12}})}^{p}\leq
C\sum\limits_{s=-4}^{\infty}\sum\limits_{j=0}^{\infty}
\sum\limits_{s_0=2}^{s+j+6}
\sum\limits_{Q_k\in\mathcal{F}^{s_0}}
2^{(j+|m|)\tau p -j\beta p-\frac{2mp}{np'}}
||f||^p_{L^p(\mathcal{F}^{s,j}_{Q_k})}.$$
Since, by Lemma \ref{g}, for any fixed $s,j$ and $s_0$,
$$\sum\limits_{Q_k\in\mathcal{F}^{s_0}}
||f||^p_{L^p(\mathcal{F}^{s,j}_{Q_k})}\leq
\sum\limits_{Q_k\in\mathcal{F}^{s_0}}\sum\limits_{Q_l\in\mathcal{F}^{s,j}_{Q_k}}
||f||^p_{L^p({Q_l})}
\leq C2^{j(n-1)}||f||^p_{L^p(\mathcal{F}^{s})},$$
we have
\begin{equation}\label{add}
\begin{array}{ccc}
||D^2 w||_{L^p(\Omega_{1/{12}})}^{p}
\leq
C\displaystyle\sum\limits_{s=-4}^{\infty}\sum\limits_{j=0}^{\infty}
\sum\limits_{s_0=2}^{s+j+6}
2^{(j+|m|)\tau p-j\beta p-\frac{2mp}{np'}+j(n-1)}
||f||^p_{L^p(\mathcal{F}^{s})}\\
= C\displaystyle\sum\limits_{s=-4}^{\infty}
\left(||f||^p_{L^p(\mathcal{F}^{s})}\cdot
\sum\limits_{j=0}^{\infty}
\sum\limits_{s_0=2}^{s+j+6}
2^{(j+|m|)\tau p-j\beta p-\frac{2mp}{np'}+j(n-1)}
\right).
\end{array}
\end{equation}

Now we choose $\tau$ as follows.
Recall that $\beta=\alpha+\frac{n}{p}+\frac{2}{np'}-1$. Since $\alpha>1-1/p$, we have $\beta-\frac{2}{np'}=\alpha+\frac{n}{p}-1>\frac{n-1}{p}$ and then take $\tau>0$ small enough such that
$$\beta-\frac{2}{np'}-2\tau>\frac{n-1}{p}\ \ \mathrm{and} \ \ \frac{2}{np'}>\tau.$$
Set $$\sigma=\min \{p(\beta-\frac{2}{np'}-2\tau)-(n-1),\ p(\frac{2}{np'}-\tau)\}>0.$$

As $2\leq s_0\leq s$, we have $|m|=s-s_0=m$ and then
$$2^{(j+|m|)\tau p-j\beta p-\frac{2mp}{np'}+j(n-1)}=2^{-j(p(\beta-\tau)-(n-1))-mp(\frac{2}{np'}-\tau)}\leq 2^{-(j+m)\sigma}.$$

As $s+1\leq s_0\leq s+j+6$, we have $|m|=s_0-s=-m\leq j+6$ and then
$$2^{(j+|m|)\tau p-j\beta p-\frac{2mp}{np'}+j(n-1)}\leq
2^{(2j+6)\tau p-j\beta p+\frac{2(j+6)p}{np'}+j(n-1)}\leq 2^{-j\sigma+6\tau p+\frac{12 p}{np'}}\leq C2^{-j\sigma},$$
where $C$ depends on $n$ and $p$.

From the above two estimates, it follows that,
$$\begin{array}{c}\displaystyle\sum\limits_{j=0}^{\infty}
\sum\limits_{s_0=2}^{s+j+6}
2^{(j+|m|)\tau p-j\beta p-\frac{2mp}{np'}+j(n-1)}\leq
\sum\limits_{j=0}^{\infty}
\left(\sum\limits_{s_0=2}^{s}
2^{-(j+m)\sigma}+
C\sum\limits_{s_0=s+1}^{s+j+6}
2^{-j\sigma}\right)\\
\displaystyle=
\sum\limits_{j=0}^{\infty}
\left(\sum\limits_{s_0=2}^{s}
2^{-(j+s-s_0)\sigma}+
C(j+6)
2^{-j\sigma}\right)
\leq C.
\end{array}$$

Substitute it into \eqref{add} and then
$$
||D^2 w||_{L^p(\Omega_{1/{12}})}^{p}\leq C
\displaystyle\sum\limits_{s=-4}^{\infty}
||f||^p_{L^p(\mathcal{F}^{s})}\leq C ||f||^p_{L^p(\Omega_1)},$$
where $C$ depends on $n,\alpha,p$ and $||\varphi||_{C^{1,\alpha}(B_1')}$.
Thus, \eqref{w} holds. Combining \eqref{v} and \eqref{w}, we conclude (4.2).
\end{proof}

\section{Fully nonlinear elliptic equation}
In this section, we will exploit our method to fully nonlinear elliptic equations and
the main result is the following theorem.

\bigskip

\begin{thm}\label{t5.1}
Let $1<p<\infty$ and $0<\alpha_0\leq\alpha\leq 1$. Assume that $\Omega$ is of class $C^{1,\alpha}$ with $0\in \partial\Omega$ and $u$ is a solution of the following elliptic equation
\begin{equation}\label{FF}	
F(D^2 u,x)=f(x)\ \ \mathrm{in}\ \ \Omega_1 \ \ \mathrm{with}\ \ f\in L^{p}(\Omega_1).
\end{equation}
Suppose $F$ satisfies interior $W^{2,p}$ estimate with constant $c_e$,  that is, for any solution $v$ of \eqref{FF} and
any $B_r(x_0)\subset\Omega_1$,
\begin{equation}\label{ia}
||D^2v ||_{ L^{p}(B_{r/2}(x_0))}\leq c_e\left( r^{{n}/p-2}||v||_{L^{\infty}(B_r(x_0))}+ ||f||_{L^{p}(B_r(x_0))}\right)
\end{equation}
and $u$ satisfies pointwise boundary $C^{1,\alpha_0}$ estimate with constant $c_b$, that is, for any $x_0\in (\partial\Omega)_{1/2}$, there exists an affine function $l_{x_0}$ such that
\begin{equation}\label{ba}
|(u-l_{x_0})(x)|\leq c_b|x-x_0|^{1+\alpha_0}\ \ \mathrm{and}\ \ |Dl_{x_0}|\leq c_b.
\end{equation}
Then we have the following two estimates:

$(i)$ If $\alpha_0>1-1/p$, then
\begin{equation}\label{5.11}
||D^2 u ||_{L^{p}(\Omega_{1/{12}})}\leq C\left(1+ ||f||_{L^{p}(\Omega_1)}\right),
\end{equation}
where $C$ depends on $n$, $\alpha_0$, $p$, $c_e$, $c_b$ and $\Omega$.

$(ii)$ If $\alpha_0\leq1-1/p$, then for any $1\leq p_0<1/(1-\alpha_0)$,
\begin{equation}\label{5.120}
||D^2 u ||_{L^{p_0}(\Omega_{1/{12}})}\leq C\left(1+ ||f||_{L^{p}(\Omega_1)}\right),
\end{equation}
where $C$ depends on $n$, $\alpha_0$, $p$, $p_0$, $c_e$, $c_b$ and $\Omega$.
\end{thm}

\begin{proof}
Let $\{Q_k\}_{k=1}^{\infty}$ be Whitney decomposition of $\Omega_{1}$, ${\widetilde Q_k}=\frac65 Q_k$ and we first
prove that for any  $Q_k\subset\widetilde Q_k\subset \Omega_{1/4}$,
\begin{equation}\label{5.12}
\begin{array}{ccc}
||D^2 u||_{L^p(Q_k)}\leq  C (d_k^{n/p+\alpha_0-1}
+||f||_{L^{p}(\widetilde Q_k)}),
\end{array}
\end{equation}
where $C$ depends on $c_e$ and  $c_b$.

Indeed, since $\widetilde Q_k\subset \Omega_{1/4}$, there exist two points $y_k\in (\partial\Omega)_{1/2}$ and $\tilde x_k\in \partial\widetilde Q_k$ such that
$$|\tilde x_k-y_k|=\dist(\widetilde Q_k, \partial\Omega_1)<\dist(Q_k, \partial\Omega_1)\leq 4d_k,$$
where Lemma \ref{l2.1} (iii) is used in the last inequality.
Consequently,
$$ |x-y_k|\leq|x-\tilde x_k|+|\tilde x_k-y_k|\leq 6d_k,\ \ \forall x\in\widetilde Q_k.$$
By \eqref{ba}, there exists an affine function $l_{y_k}$(written as $l$ for simplicity in the following) such that
\begin{equation*}
|(u-l)(x)|\leq  c_b |x-y_k|^{1+\alpha_0}\leq c_b(6d_k)^{1+\alpha_0},
\ \forall x\in \widetilde Q_k.
\end{equation*}
Since $u-l$ still satisfies $F(D^2(u-l),x)=f(x)$,
we have, by \eqref{ia},
\begin{equation*}
\begin{aligned}
||D^2 (u-l)||_{L^p(Q_k)}
&\leq c_e (d_k^{n/p-2}||u-l||_{L^{\infty}(\widetilde Q_k)}+||f||_{L^{p}(\widetilde Q_k)})\\
&\leq C (d_k^{n/p+\alpha_0-1}
+||f||_{L^{p}(\widetilde Q_k)}),
\end{aligned}
\end{equation*}
where $C$ depends on $c_e$ and $c_b$. Thus, \eqref{5.12} holds.

For any $1\leq q\leq p$, by H$\rm{\ddot{o}}$lder's inequality, \eqref{5.12} and Young's inequality, we deduce
\begin{equation}\label{5.4}
\begin{aligned}
\displaystyle\int_{Q_k}|D^2 u|^{q}dx
&\leq Cd_k^{n(1-q/p)}
\left(\displaystyle\int_{Q_k}|D^2 u|^{p}dx\right)^{q/p}\\
&\leq C\left\{d_k^{n-(1-\alpha_0){q}}+
d_k^{n-{q} n/p}
\left(\displaystyle\int_{\widetilde  Q_k}|f|^{p}dx\right)^{{q}/p}\right\}\\
&\leq C\left(d_k^{n-(1-\alpha_0){q}}+
d_k^n+
\displaystyle\int_{\widetilde  Q_k}|f|^{p}dx\right),
\end{aligned}
\end{equation}
where $C$ depends on $n$, $p$, $q$, $c_e$ and  $c_b$.

Now we are ready to show \eqref{5.11} and \eqref{5.120}.

If $\alpha_0> 1-1/p$, then $n-(1-\alpha_0){{p}}>n-1$ and by Lemma \ref{lf},
$$\displaystyle\sum\limits_{\widetilde Q_k\subset \Omega_{1/4}}d_k^{n-(1-\alpha_0){p}}\leq C,$$
where $C$ depends on $n,\alpha_0,p$ and $\Omega$.
Set $q=p$ in \eqref{5.4} and it follows that
$$\displaystyle\int_{Q_k}|D^2 u|^{p}dx
\leq C\left(d_k^{n-(1-\alpha_0){p}}+
d_k^n+
\displaystyle\int_{\widetilde  Q_k}|f|^{p}dx\right).$$
Since $\Omega_{1/{12}}\subset\bigcup_{\widetilde Q_k\subset \Omega_{1/4}}Q_k$, we deduce from the above two estimates that
\begin{equation*}
\begin{aligned}
&\displaystyle\int_{\Omega_{1/{12}}}|D^2 u|^{{p}}dx
\leq \displaystyle\sum\limits_{\widetilde Q_k\subset \Omega_{1/4}}
\displaystyle\int_{Q_k}|D^2 u|^{{p}}dx\\
&\leq C\displaystyle\sum\limits_{\widetilde Q_k\subset \Omega_{1/4}}
\left(d_k^{n-(1-\alpha_0){p}}+
d_k^n+
\displaystyle\int_{\widetilde  Q_k}|f|^{p}dx\right)
\leq C\left(1+\displaystyle\int_{\Omega_1}|f|^{p}dx\right),
\end{aligned}
\end{equation*}
where $C$ depends on $n$, $\alpha_0$, $p$, $c_e$, $c_b$ and $\Omega$. This gives \eqref{5.11}.

If $\alpha_0\leq1-1/p$, then $p\geq 1/(1-\alpha_0)$.
For any $1\leq p_0<1/(1-\alpha_0)$, we have $n-(1-\alpha_0){{p_0}}>n-1$ and then by Lemma \ref{lf},
$$\displaystyle\sum\limits_{\widetilde Q_k\subset \Omega_{1/4}}d_k^{n-(1-\alpha_0){p_0}}\leq C,$$
where $C$ depends on $n,\alpha_0,p_0$ and $\Omega$.
Since now $p_0<p$,  set $q=p_0$ in \eqref{5.4} and we obtain
$$\displaystyle\int_{Q_k}|D^2 u|^{p_0}dx
\leq C\left(d_k^{n-(1-\alpha_0){p_0}}+
d_k^n+
\displaystyle\int_{\widetilde  Q_k}|f|^{p_0}dx\right).$$
From $\Omega_{1/{12}}\subset\bigcup_{\widetilde Q_k\subset \Omega_{1/4}}Q_k$ and the above two estimates, we deduce that
\begin{equation*}
\begin{aligned}
&\displaystyle\int_{\Omega_{1/{12}}}|D^2 u|^{{p_0}}dx
\leq \displaystyle\sum\limits_{\widetilde Q_k\subset \Omega_{1/4}}
\displaystyle\int_{Q_k}|D^2 u|^{{p_0}}dx\\
&\leq C\displaystyle\sum\limits_{\widetilde Q_k\subset \Omega_{1/4}}
\left(d_k^{n-(1-\alpha_0){p_0}}+
d_k^n+
\displaystyle\int_{\widetilde  Q_k}|f|^{p}dx\right)
\leq C\left(1+\displaystyle\int_{\Omega_1}|f|^{p}dx\right),
\end{aligned}
\end{equation*}
where $C$ depends on $n$, $\alpha_0$, $p$, $p_0$, $c_e$, $c_b$ and $\Omega$. Hence \eqref{5.120} holds.
\end{proof}

\bigskip

\begin{re}
From Theorem 5.1, we see again that  by Whitney decomposition, local boundary $W^{2,p}$ estimate follows from interior $W^{2,p}$ estimate and boundary $C^{1,\alpha}$ estimate which are assumed. As for interior $W^{2,p}$ estimate, we refer to \cite{C} and Theorem 7.1 in \cite{CC}; as for
boundary $C^{1,\alpha}$ estimate, we refer to \cite{LZ}  and \cite {SS}.


\end{re}

\bigskip

\section*{Acknowledgement}
This work is supported by NSFC 12071365.




\bigskip

\bibliographystyle{elsarticle-num}



\end{document}